\documentclass[12pt,twoside]{article}

\begin{document}

\date{}




\centerline {\Large{\bf Well-Posedness Of Common Fixed Point   }}

\centerline{}

\centerline {\Large{\bf Theorem For Three and Four Mappings }}

\centerline{}

\centerline {\Large{\bf  Under Strict Contractive Conditions }}

\centerline{}

\centerline {\Large{\bf In Fuzzy metric Space}}

\centerline{}

\centerline{\bf { Sumit Mohinta and T. K. Samanta }}
\centerline{Department of Mathematics, Uluberia College, India-711315.} %
\centerline{e-mail: mumpu$_{-}$tapas5@yahoo.co.in}
\centerline{e-mail: sumit.mohinta@yahoo.com}

\begin{abstract}
\textbf{\emph{None has studied the well-posedness of common fixed points in fuzzy metric space. In this paper, our target is to develop the well-posedness of common fixed points in fuzzy metric space. Also using weakly compatibility, implicit relation, property (E.A.) and strict contractive conditions, we have established the unique common fixed point for three self mappings and also for four self mappings in fuzzy metric space.}}
\end{abstract}


\newtheorem{Theorem}{\quad Theorem}[section]

\newtheorem{Definition}[Theorem]{\quad Definition}

\newtheorem{Corollary}[Theorem]{\quad Corollary}

\newtheorem{Lemma}[Theorem]{\quad Lemma}

\newtheorem{Note}[Theorem]{\quad Note}

\newtheorem{Remark}[Theorem]{\quad Remark}

\newtheorem{Result}[Theorem]{\quad Result}

\newtheorem{Proposition}[Theorem]{\quad Proposition}

\newtheorem{Example}[Theorem]{\quad Example}



\section{Introduction}


Motivated by a work due to Popa\cite{Popa}, different authors\cite{Aliouche,Jungck,Pathak,Akkouchi} have tried to prove fixed point theorems using an implicit relation, which is a good idea since it covers several contractive conditions rather than one contractive condition in an ordinary metric space. In fact, it is seen that commuting implies weakly commuting which also implies compatible and there are examples in the literature verifying that the inclusions are proper, see\cite{Jungck}. In the paper\cite{Jungck1}, Jungck defined the weakly compatible maps and established that two maps are weakly compatible if they commute at their coincidence points. Using the concept of weakly compatible maps, implicit relation, property(E.A.), the authors\cite{Akkouchi} established the unique common fixed point for three self mappings under strict contractive conditions\cite{Aamri} in an ordinary metric space. Also the authors\cite{Akkouchi} proved that such type fixed point problem is well posed. The author, A. Aliouche\cite{Aliouche}, also established the unique common fixed point for four self mappings using the concept of weakly compatible maps, implicit relation, property(E.A.) and strict contractive condition in an ordinary metric space, but he didn't establish the well-posedness of such type common fixed point.\\
So far our knowledge, a little bit of such type results have been developed in fuzzy metric space, a lot of common fixed point theorems for three and four self mappings and their well-posedness yet remain to develop in fuzzy metric space. \\
Fuzzy set theory was first introduce by Zadeh\cite{zadeh} in 1965 to describe the situation in which data are imprecise or vague or uncertain. It has a wide range of application in the field of population dynamics , chaos control , computer programming , medicine , etc.\\
The concept of fuzzy metric was first introduced by Kramosil and Michalek\cite{Kramosil} and later on it is modified and a few concepts of mathematical analysis have been developed by George and Veeramani\cite{Veeramani,Veeramani1} and also they have developed the fixed point theorem in fuzzy metric space\cite{Veeramani2}. In fuzzy metric space, the notion of compatible maps under the name of asymptotically commuting maps was introduced in the paper\cite{Mishra} and then in the paper\cite{Singh}, the notion of weak compatibility has been studied in fuzzy metric space. However , the study of common fixed points of non-compatible maps is of great interest, which has been initiated by Pant. With the help of the property (E.A.), which was introduced in the paper\cite{Aamri}, Pant and Pant\cite{Pant} studied the common fixed points of a pair of non-compatible maps in fuzzy metric space. \\ None has studied the well-posedness of common fixed points in fuzzy metric space. In this paper, our target is to develop the well-posedness of common fixed points in fuzzy metric space. Also using weakly compatibility, implicit relation, property(E.A.) and strict contractive conditions, we have established the unique common fixed point for three self mappings and also for four self mappings in fuzzy metric space.


\section{Preliminaries}

We quote some definitions and statements of a few theorems which will be
needed in the sequel.

\begin{Definition}\cite{Schweizer}
A binary operation \, $\ast \; : \; [\,0 \; , \; 1\,] \; \times \;
[\,0 \; , \; 1\,] \;\, \longrightarrow \;\, [\,0 \; , \; 1\,]$ \, is
continuous \, $t$ - norm if \,$\ast$\, satisfies the
following conditions \, $:$ \\
$(\,i\,)$ \hspace{0.5cm} $\ast$ \, is commutative and associative ,
\\ $(\,ii\,)$ \hspace{0.3cm} $\ast$ \, is continuous , \\
$(\,iii\,)$ \hspace{0.2cm} $a \;\ast\;1 \;\,=\;\, a \hspace{1.2cm}
\forall \;\; a \;\; \varepsilon \;\; [\,0 \;,\; 1\,]$ , \\
$(\,iv\,)$ \hspace{0.3cm} $a \;\ast\; b \;\, \leq \;\, c \;\ast\; d$
\, whenever \, $a \;\leq\; c$  ,  $b \;\leq\; d$  and  $a \, , \, b
\, , \, c \, , \, d \;\, \varepsilon \;\;[\,0 \;,\; 1\,]$.
\end{Definition}

\begin{Result}\cite{klement}
$(\,a\,)\;$  For any \, $r_{\,1} \; , \;
r_{\,2} \;\; \in\;\; (\,0 \;,\; 1\,)$ \, with \, $r_{\,1} \;>\;
r_{\,2}$, there exist $ \; r_{\,3} \;\; \in \;\; (\,0
\;,\; 1\,)$ \, such that \, $r_{\,1} \;\ast\; r_{\;3} \;>\; r_{\,2}$ ,
\\ $(\,b\,)$ \, For any \, $r_{\,5} \;\,
\in\;\, (\,0 \;,\; 1\,)$ , there exist \, $r_{\,6} \; \, \in\;\, (\,0 \;,\; 1\,)$ \, such that  ${\hspace{1.5cm}}r_{\,6} \;\ast\;
r_{\,6} \;\geq\; r_{\,5}$.
\end{Result}

\begin{Definition}\cite{Veeramani}
The $3$-tuple $(\,X \,,\, \mu \,,\, \ast\,)$ is called a fuzzy metric space if X is an arbitrary set, \,$\ast$\, is a continuous t-norm and \,$\mu$ is a fuzzy set in \,$X^{2}\,\times\,(0,\infty)$\, satisfying the following conditions \,:
\\$(\,i\,)$ \hspace{0.5cm}$\mu\,(\,x \;,\,y \;,\; t\,) \;\,>\;\, 0 \, ;$ \\
$(\,ii\,)$ \hspace{0.2cm} $\mu\,(\,x \;,\,y\;,\; t\,) \;\,=\;\, 1$ \, if
and only if \, $x \;=\;y \,$\,
\\$(\,iii\,)$ \hspace{0.1cm}$\mu\,(\,x \;,\,y\;,\; t\,) \;\,=\;\,\mu\,(\,y \;,\,x\;,\; t\,);$
\\ $(\,iv\,)$\hspace{0.1cm} $\mu\,(\,x \;,\,y \;,\; s\,) \;\ast\; \mu\,(\,y \;,\,z \;,\; t\,)
\;\,\leq\;\, \mu\,(\,x \;,\,z \;,\; s\;+\;t\,
\,) \, ;$
\\$(\,v\,)$ \hspace{0.2cm}$\mu\,(\,x \,,\,y\,,\,\cdot\,) :(0 \,,\;\infty\,)\,\rightarrow
\,(0 \,,\;1]$ \, is continuous for all x\,,\,y\,,\,z $\,\in\,X$ and ${\hspace{1.4cm}}t,\,s\,>\,0.$
\end{Definition}

\begin{Definition}\cite{Samanta}
Let \,$(\,X \,,\, \mu \,,\, \ast\,)$\, be a fuzzy metric space.
A sequence $\{\,x_{\,n}\,\}_{\,n}$ in X is called Cauchy sequence if and only if
\[\mathop {\lim }\limits_{n\,\, \to \,\,\infty } \,\mu\,(\,x_{n} \,,\,x_{n\,+\,p} \,,\, t\,)\;=\;1\;for \,\,each\,\, t \,>\, 0 \;and\; p \,=\, 1 \,,\, 2 \,,\, 3 \,,\, \cdots\]
A sequence \,$\{\,x_{\,n}\,\}_{\,n}$\, in \,$X$\, is said to converge to
\,$x\,\in\,X$\, if and only if
\[\mathop {\lim }\limits_{n\,\, \to \,\,\infty } \,\mu\,(\,x_{n} \,,\,x \,,\, t\,)\;=\;1 \;\;for \,\,each\,\, t \,>\, 0 \hspace{4.5cm}\]
A fuzzy metric space \,$(\,X \,,\, \mu \,,\, \ast\,)$\, is said to be complete if and only if every Cauchy sequence in \,$X$\, is convergent in \,$X$ .
\end{Definition}

\begin{Definition}\cite{Mishra}
 Let $A$ and $B$ be maps from an fuzzy metric space \,$(\,X \,,\, \mu \,,\, \ast\,)$\, into itself. The maps $A$ and $B$ are said to be compatible $($ or asymptotically commuting $)$, if for all \,$t \,>\, 0$ ,
 \[\mathop {\lim }\limits_{n\,\, \to \,\,\infty } \,\mu\,(\,A\,Bx_{n} \,,\,B\,Ax_{n\,} \,,\, t\,)\;=\;1\]
whenever \,$\{\,x_{\,n}\,\}$\, is a sequence in \,$X$\, such that $\mathop {\lim }\limits_{n\,\, \to \,\,\infty }\,Ax_{\,n}\,=\,\mathop {\lim }\limits_{n\,\, \to \,\,\infty }Bx_{\,n} \,=\, z$
for some \,$z\,\in\,X.$

\end{Definition}

\begin{Definition}\cite{Singh}
Let $A$ and $B$ be maps from a fuzzy metric space \,$(\,X \,,\, \mu \,,\,\\ \ast\,)$\, into itself. The maps are said to be weakly compatible if they commute at their coincidence points, that is, $Az\,= \,Bz$\, implies that \,$ABz\,=\,BAz.$
\\Note that compatible mappings are weakly compatible but converse is not true in general.
\end{Definition}

\begin{Definition}\cite{Pant}
Let $A$ and $B$ be two self-maps of a fuzzy metric space \,$(\,X \,,\, \mu \,,\, \ast\,)$\,. We say that $A$ and $B$ satisfy the property $(\,E.A.\,)$ if there exists a sequence \,$\{x_{\,n}\}$\, such that
\[\mathop {\lim }\limits_{n\,\, \to \,\,\infty }\,Ax_{n}\,=\,\mathop {\lim }\limits_{n\,\, \to \,\,\infty }Bx_{n}\,=\,z\]
for some \,$z \,\in\, X$. Note that weakly compatible and property $(\,E.A.\,)$ are independent to each other $(\,See\cite{Pathak}, Ex.2.2\,)$
\end{Definition}

\begin{Definition}\cite{Pant}
The mapping $A \,,\, B \,,\, S \,,\, T \,:\, \longrightarrow\, X$ of a fuzzy metric space \,$(\,X \,,\, \mu \,,\, \ast\,)$\, satisfy a common property $(\,E.A.\,)$ if there exist two sequences \,$\{\,x_{\,n}\,\}$\, and \,$\{\,y_{\,n}\,\}$\,  such that
\[\mathop {\lim }\limits_{n\,\, \to \,\,\infty }\,Ax_{\,n}\,=\,\mathop {\lim }\limits_{n\,\, \to \,\,\infty }Sx_{\,n}\,=\,\mathop {\lim }\limits_{n\,\, \to \,\,\infty }\,By_{\,n}\,=\,\mathop {\lim }\limits_{n\,\, \to \,\,\infty }Ty_{\,n}\,=\,z\]
for some \,$z \,\in\, X$. If \,$B\,=\,A $\, and \,$T\,=\,S$, are obtain definition \,$(\,2.7\,)$\,
\end{Definition}

\begin{Definition}\cite{Samanta}
Let \,$(\,X \,,\, \mu \,,\, \ast\,)$\, be a fuzzy metric space. A subset \,$P\, of\, X$\, is said to be closed if for any sequence \,$\{\,x_{\,n}\,\}$\, in \,$P$\, converges to \,$x\,\in\,P,$\, that is
\[\mathop {\lim }\limits_{n\,\, \to \,\,\infty } \,\mu\,(\,x_{\,n} \,,\,x \,,\, t\,)\;=\;1{\hspace{0.5cm}}\Longrightarrow\;\,x\,\in\,P{\hspace{1.0cm}}\forall\;t\,>\,0\]
\end{Definition}

\section{ Implicit Relations}


\begin{Definition}
Let \,$I\, = \,[\,0 \,,\, 1\,]$\, and \,$F:\,I^{\,6}\,\rightarrow\,I$\, be continuous function. We define the following property:
\[(\,F_{\,1}\,)\,:\, F\,(\,u(\,t\,)\,,\,1\,,\,1\,,\,u(\,t\,)\,,\,u(\,t\,)\,,\,1\,)\,>\,1{\hspace{0.8cm}}
\forall\;t\,>\,0 \,,\, 0\,\leq\,u(\,t\,)\,<\,1\]
\[(\,F_{\,2}\,):\, F\,(\,u(\,t\,)\,,\,1\,,\,u(\,t\,)\,,\,1\,,\,1\,,\,u(\,t\,)\,)\,>\,1{\hspace{0.8cm}}
\forall\,t\,>\,0 \,,\, 0\,\leq\,u(\,t\,)\,<\,1\]
\[(\,F_{\,3}\,):\, F\,(\,u(\,t\,)\,,\,u(\,t\,)\,,\,1\,,\,1\,,\,u(\,t\,)\,,\,u(\,t\,)\,)\,>\,1
{\hspace{0.3cm}}\forall\,t\,>\,0 \,,\, 0\,\leq\,u(\,t\,)\,<\,1\]
\end{Definition}

\begin{Example}Let \,$F\,(t_{\,1} \,,\, \cdots \,,\, t_{\,6})\,:\, \frac{t_{\,1}\,+\,t_{\,2}\,+\,t_{\,3}}{k\;max\,{\{t_{\,4} \,,\, t_{\,5} \,,\, t_{\,6}\}}}$\, where \,$k\,\in\,(\,0 \,,\, 1\,)$

\[(\,F_{\,1}\,)\,:\,F\,(\,u(\,t\,)\,,\,1\,,\,1\,,\,u(\,t\,)\,,\,u(\,t\,)\,,\,1\,)
{\hspace{5.0cm}}\]
\[{\hspace{3.5cm}}=\,\frac{u(\,t\,)\,+\,1\,+\,1}{k\,max\,{\{u(\,t\,),\,
u(\,t\,),\,1\}}}\,>\,
1{\hspace{1.0cm}}\forall\,t\,>\,0,\,0\,\leq\,u(\,t\,)\,<\,1\]
\[(\,F_{\,2}\,)\,:\,F\,(\,u(\,t\,),\,1,\,u(\,t\,),\,1,\,1,\,
u(\,t\,)\,){\hspace{5.4cm}}\]
\[{\hspace{3.1cm}}=\,\frac{u(\,t\,)\,+\,1\,+\,
u(\,t\,)}{k\,max{\{\,1,\,1,\,u(\,t\,)\}}}\,>\,
1{\hspace{1.0cm}}\forall\,t\,>\,0,\,0\,\leq\,u(\,t\,)\,<\,1\]
\[(\,F_{\,3}\,)\,:\,F\,(\,u(\,t\,),\,u(\,t\,),\,1,\,1,\,
u\,(t\,),\,u\,(\,t\,)\,){\hspace{5.0cm}}\]
\[{\hspace{3.0cm}}=\,\frac{u(\,t\,)\,+\,
u(\,t\,)\,+\,1}{k\,max\,{\{\,1,\,u(\,t\,),\,u(\,t\,)\}}}\,>\,
1{\hspace{0.6cm}}\forall\,t\,>\,0,\,0\,\leq\,u(\,t\,)\,<1\]
\end{Example}

\begin{Definition}
Let \,$I \,=\, [\,0 \,,\, 1\,]$\, and \,$F \,:\, I^{\,6}\, \rightarrow\, I$\, be continuous function. We define the following property\,:
\\$(\,F_{\,4}\,) $\, There exists \,$k \,\in\, (\,0 \,,\, 1\,)$\, such that \,$0 \,\leq\, u(\,t\,) \,,\, v(\,t\,) \,,\, w(\,t\,) \,<\, 1 \,,$ \[F\,(\,u(\,t\,) \,,\, v(\,t\,) \,,\, 1 \,,\, w(\,t\,) \,,\, u(\,t\,) \,,\,
v(\,t\,)\,) \,\leq\, 1{\hspace{1.0cm}}\forall\;t \,>\, 0 ,\]
\\${\hspace{2.5cm}}\Longrightarrow\; u(\,t\,) \,\geq\, \frac{1}{k}\;w(\,t\,)$
\end{Definition}

\begin{Example}Let $F\,(t_{\,1} \,,\, \cdots \,,\, t_{\,6}) \,:\, max\,{\{t_{\,2} \,,\, t_{\,3} \,,\, t_{\,5} \,,\, t_{\,6}\}} \,-\, k\,t_{\,1} \,+\, t_{\,4}$
\\ where \,$k \,\in\, (\,0 \,,\, 1\,)$
\[F\,(\,u(\,t\,) \,,\, v(\,t\,) \,,\, 1 \,,\, w(\,t\,) \,,\, u(\,t\,) \,,\,
v(\,t\,)\,) \,\leq\, 1{\hspace{1.0cm}}\forall\;t \,>\, 0,\]
\[\Longrightarrow\; max\,\{v(\,t\,) \,,\, 1 \,,\, u(\,t\,) \,,\,
v(\,t\,)\,\} \,-\, k\,u(\,t\,) \,+\, w(\,t\,) \,\leq\, 1\]
\[\Longrightarrow\; 1 \,-\, k\,u(\,t\,) \,+\, w(\,t\,) \,\leq\, 1{\hspace{4.9cm}}\]
\[\Longrightarrow\; u(\,t\,) \,\geq\, \frac{1}{k}\;w(\,t\,){\hspace{6.5cm}}\]
\end{Example}

\begin{Example}Let F $(t_{\,1} \,,\, \cdots \,,\, t_{\,6}):\, \;t_{1}\,+\,\,t_{2} \,+\, \,t_{3} \,+\, \frac{1}{k}\,t_{4} \,-\, 2\,t_{5} \,-\, t_{6}$
\\ where \,$k \,\in\, (\,0 \,,\, 1\,)$
\[F\,(\,u(\,t\,) \,,\, v(\,t\,) \,,\, 1 \,,\, w(\,t\,) \,,\, u(\,t\,) \,,\,
v(\,t\,)\,) \,\leq\, 1{\hspace{1.0cm}}\forall\;t \,>\, 0,\]
\[\Longrightarrow\;u(t)\,+\,v(t)\,+\,1\,+\,\frac{1}{k}\,w(t)\,-\,2\,u(t)\,-\,v(t)\,\leq\,1{\hspace{0.6cm}}\]
\[\Longrightarrow\;u(t)\,\geq\,\frac{1}{k}\,w(t){\hspace{6.7cm}}\]
\end{Example}


\section{Common Fixed Point}

\begin{Theorem}
Let \,$A \,,\, B$ \,and\, $I$\, be three self mappings of a fuzzy metric space \,$(\,X \,,\, \mu \,,\, \ast\,)$\, such that:
\\${\hspace{2.0cm}}(\,i\,)$ The pair \,$\{\,A \,,\, I\,\}$\, and \,$\{\,B \,,\, I\,\}$\, are weakly compatible,
\\${\hspace{2.0cm}}(\,ii\,)$ The mapping \,$A \,,\, B$ \,and\, $I$\, satisfy the property $(\;E.A.\;),$
\\${\hspace{2.0cm}}(\,iii\,) \,F\,(\,\mu(\,Ax \,,\, By \,,\, t\,) \,,\, \mu(\,Ix \,,\, Iy\, ,\, t\,) \,,\, \mu(\,Ix \,,\, Ax \,,\, t\,) \,,\, \mu(\,Iy \,,\, By \,,\, t\,),
\\{\hspace{5.0cm}}\mu(\,Ix \,,\, By \,,\, t\,) \,,\, \mu(\,Iy \,,\, Ax \,,\, t\,)\,) \,\leq\, 1\,\,\forall\,x\neq\,y\,\;in\,X$,
\\Where F satisfies property $(\,F_{\,1}\,) \,,\, (\,F_{\,2}\,)$ \,and\, $(\,F_{\,3}\,)$,
\\${\hspace{2.0cm}}(\,iv\,)$ \,$I(\,x\,)$\,is closed,
\\\\Then the mappings \,$A \,,\, B$ \,and\, $I$\, have a unique common fixed point.
\end{Theorem}

{\bf Proof.} Since the set of mappings \,$\{\,A\,,\,B\,,\,I\,\}$\, satisfies the property \,$(\;E.A\;),$\, there exists a sequence \,$\{\,x_{\,n}\}$\, such that
\[\mathop {\lim }\limits_{n\,\, \to \,\,\infty }\,Ax_{\,n}\,=\,\mathop {\lim }\limits_{n\,\, \to \,\,\infty }Bx_{\,n}\,=\,\mathop {\lim }\limits_{n\,\, \to \,\,\infty }\,Ix_{\,n}\,=\,u{\hspace{1.0cm}}(1)\]
for some \,$u\in\,X.$
\\Since \,$I(\,X\,)$\, is closed there exists a point \,$a\in X$\, such that \,$u\,=\,Ia$\,. If the sequence \,$\{\,x_{\,n}\,\}$ satisfies
\[x_{\,n}\,=\,a\,,\forall\;n\;\geq n_{\,0}\]
for some positive integer \,$n_{\,0}$ , then from \,$(\,1\,)$\,, we have
\[u\,=\,Aa\,=\,Ia\,=\,Ba\]
So , we may suppose that \,$x_{\,n}\,\neq \,a$\, for all integer n , (otherwise , we consider a subsequence satisfying this property). By putting \,$x\,=\,x_{\,n}$\, and \,$y\,=\,a$\, in \,$(\,iii\,)$\, we obtain:
\[\,F(\,\mu(\,Ax_{\,n}\,,\,Ba\,,\,t\,)\,,\,\mu(\,Ix_{\,n}\,,\,Ia\,,\,t\,)
\,,\,\mu(\,Ix_{\,n}\,,\,Ax_{\,n}\,,\,t\,)\,,\,
\mu(\,Ia\,,\,Ba\,,\,t\,)\,,\]
\[{\hspace{4.0cm}}\mu(\,Ix_{\,n}\,,\,Ba\,,\,t\,)\,,\,\mu(\,Ia\,,\,Ax_{\,n}\,,\,t\,)\,)\,
\leq\,1\,{\hspace{0.5cm}}\,\forall \;t \,>\, 0.\]
Letting \,$n\,\rightarrow\,\infty\,, $\, we obtain :
\[\,F(\,\mu(\,Ia\,,\,Ba\,,\,t\,)\,,\,1\,,1\,,\,\mu(\,Ia\,,\,Ba\,,\,t\,)
\,,\,\mu(\,Ia\,,\,Ba\,,\,t\,)\,,\,1\,,\,)\,\leq\,1\]
which, by virtue of \,$(F_{\,1})$ \,,\, implies that
\[\mu(\,Ia\,,\,Ba\,,\,t\,)\,=\,1{\hspace{0.5cm}} \forall\; t \,>\, 0.\]
\[\Rightarrow \,Ia\,=\,Ba{\hspace{3.9cm}}\]
Since \,$x_{\,n}\,\neq \,a$\, for all integers \,$n$\, and putting \,$x\,=\,a\,,\,y\,=\,x_{\,n}$\, in \,$(\,iii\,)$ ,\, then we get
\[\,F(\,\mu(\,Aa\,,\,Bx_{\,n}\,,\,t\,)\,,\,\mu(\,Ia\,,\,Ix_{\,n}\,,\,t\,)
\,,\,\mu(\,Ia\,,\,Aa\,,\,t\,)\,,\,
\mu(\,Ix_{\,n}\,,\,Bx_{\,n}\,,\,t\,)\,,\]
\[{\hspace{5.0cm}}\mu(\,Ia\,,\,Bx_{\,n}\,,\,t\,)\,,\,\mu(\,Ix_{\,n}
\,,\,Aa\,,\,t\,)\,)\,\leq\,1\,\]
Letting \,$n\,\rightarrow\,\infty\,, $\, we obtain :
\[\,F(\mu(\,Aa\,,\,Ia\,,\,t\,)\,,\,1\,,\,\mu(\,Ia\,,\,Aa\,,\,t\,)
\,,\,1\,,\,1\,,\,\mu(\,Ia\,,\,Aa\,,\,t\,)\,)
\leq\,1{\hspace{0.5cm}}\,\forall\; t \,>\, 0.\]
which, by virtue of \,$(\,F_{\,2}\,)\,,$\, implies that
$\mu(\,Aa\,,\,Ia\,,\,t\,)\,=\,1{\hspace{0.5cm}} \forall\; t \,>\, 0$
\\Hence \,$Aa \,=\, Ia$. Therefore , we obtain
\[Aa \,=\, Ia \,=\, Ba {\hspace{3.0cm}}\]
${\hspace{3.0cm}}$ we set ,
\,$x\,=\,Aa\,=\,Ia\,=\,Ba\,.{\hspace{3.0cm}}$\,
\\we shall prove that x is a common fixed point of the mappings \,$A\,,\,B$\, and\, $I.$
\\Since the pairs \,$\{\,A\,,\,I\,\}$\, and \,$\{\,B\,,\,I\,\}$\, are weakly compatible , then we have
\[AIa \,=\, IAa \,\;and \;\, BIa \,=\, IBa \hspace{1.7cm}\]
Therefore ,
\[IIa \,=\, IAa \,=\, AIa \,=\, AAa \hspace{1.7cm}\]
\[{\hspace{2.0cm}}i.e,\,IIa \,=\, AAa \,\,\Longrightarrow\,\, Ix \,=\, Ax{\hspace{1.5cm}}\cdots\hspace{1.5cm}(\,2\,)\]
and
\[{\hspace{2.0cm}} IIa \,=\, IBa \,=\, BIa \,=\, BBa \,=\, Ix \,=\, Bx{\hspace{1.5cm}}\cdots\hspace{1.5cm}(\,3\,)\]
If \,$x\,=\,a\,,$\, then we have \,$x \,=\, Ax \,=\, Ix \,=\, Bx\,.$\, Therefore x is a common fixed point of the mappings \,$A \,,\, B$\, and\, $I.$ So , we may suppose that \,$x \,\neq \, a.$ In this case , by using the equalities \,$(\,2\,)$\, and \,$(\,3\,)$\, and the inequalities \,$(\,iii\,)$\, we put \,$x\,=\,a$\, and \,$y\,=\,x$\, then
\[\,F(\,\mu(\,Aa\,,\,Bx\,,\,t\,)\,,\,\mu(\,Ia\,,\,Ix\,,\,t\,)
\,,\,\mu(\,Ia\,,\,Aa\,,\,t\,)\,,\,\mu(\,Ix\,,\,Bx\,,\,t\,)\,,\]
\[{\hspace{5.0cm}}\mu(\,Ia\,,\,Bx\,,\,t\,)\,,\,\mu(\,Ix\,,\,Aa\,,\,t\,)\,)
\,\leq\,1\,{\hspace{0.5cm}}
\forall\; t \,>\, 0\]
\[\Longrightarrow\;\; \mu(\,x\,,\,Ix\,,\,t\,)\,,\,
\mu(\,x\,,\,Ix\,,\,t\,)\,,\,1\,,\,1\,,\,\mu(\,x\,,\,Ix\,,\,t\,)\,,\,
\mu(\,x \,,\, Ix \,,\, t\,)) \,\leq\, 1 \]
a contradiction of \,$(\,F_{\,3}\,),$\, so, we have \,$\mu(\,x\,,\,Ix\,,\,t\,)\,=\,1$
\\Hence \,$Ix\,=\, x.$\, Then we get from \,$(\,2\,)$\, and \,$(\,3\,)$ \,$x\,=\,Ix\,=\,Ax\,=\,Bx\,.$
\\Therefore x is a commone fixed point of \,$A\,,\,B$\, and\, $I.$
\\Now , we show that the point x is unique common fixed point of \,$A\,,\,B$\, and\, $I.$
\\Suppose that \,$A\,,\,B$\, and\, $I.$\, have another common fixed point \,$x_{\,1}\,.$
\\Then , we put \,$y\,=\,x_{\,1}$\, in \,$(\,iii\,)$\,
\[\,F(\mu(\,Ax\,,\,Bx_{\,1}\,,\,t\,)\,,\,\mu(\,Ix\,,\,Ix_{\,1}\,,\,t\,)
\,,\,\mu(\,Ix\,,\,Ax\,,\,t\,)\,,\,
\mu(\,Ix_{\,1}\,,\,Bx_{\,1}\,,\,t\,),\]
\[{\hspace{5.0cm}}\mu(\,Ix\,,\,Bx_{\,1}\,,\,t\,)\,,\,\mu(\,Ix_{\,1}
\,,\,Ax\,,\,t\,)\,)\,\leq\,1\]
\[\Rightarrow\,F(\mu(\,x\,,\,x_{\,1}\,,\,t\,)\,,\,\mu(\,x\,,\,x_{\,1}\,,\,t\,)
\,,\,\mu(\,x\,,\,x\,,\,t\,)\,,\,
\mu(\,x_{\,1}\,,\,x_{\,1}\,,\,t\,)\,,{\hspace{1.7cm}}\]
\[{\hspace{5.0cm}}\mu(\,x\,,\,x_{\,1}\,,\,t\,)\,,\,
\mu(\,x_{\,1}\,,\,x\,,\,t\,)\,)\,\leq\,1\]
\[\Rightarrow\,F(\,\mu(\,x\,,\,x_{\,1}\,,\,t\,)
\,,\mu(\,x\,,\,x_{\,1}\,,\,t\,)\,,\,1\,,\,1\,,\,
\mu(\,x\,,\,x_{\,1}\,,\,t\,)\,,\mu(\,x_{\,1}\,,\,x\,,\,t\,)\,)\,\leq\,1\]
a contradiction of \,$(\,F_{\,3}\,),$\, we get \,$x\,=\,x_{\,1}.$\,
\\This completes the proof.

\begin{Theorem}Let \,$A\,,\,B\,,\,S$\, and \,$T$\, be self-mappings of a fuzzy metric space
\,$(\,X \,,\, \mu \,,\, \ast\,)$\, satisfying the following conditions :
\[A(\,X\,)\,\subset\,T(\,X\,) \,\;and\,\; B(\,X\,)\,\subset\,S(\,X\,)\]
$ F\,(\,\mu(\,Ax\,,\,By\,,\,t\,)\,,\,\mu(\,Sx\,,\,Ty\,,\,t\,)
\,,\,\mu(\,Ax\,,\,Sx\,,\,t\,)\,,\,
\mu(\,By\,,\,Ty\,,\,t\,)\,,\,
\\{\hspace{5.5cm}}\mu(\,Sx\,,\,By\,,\,t\,)\,,\,
\mu(\,Ax\,,\,Ty\,,\,t\,)\,) \;\leq\; 1 {\hspace{1.5cm}}\cdots{\hspace{0.5cm}}(\,4\,)$
\\\\for all \,$\,x\,,\,y$ \,in\, $X$\, and where \,$F$\, satisfies property $(\,F_{\,1}\,)\,,\,(\,F_{\,2}\,)$ \,and\, $(\,F_{\,3}\,).$
\\Suppose that \,$(\,A\,,\,S\,)$\, or \,$(\,B\,,\,T\,)$\, satisfies property \,$(\,E.A.\,)$\, and the pairs \,$(\,A\,,\,S\,)$\,and \,$(\,B\,,\,T\,)$ are weakly compatible . If the range of one \,$A\,,\,B\,,\,S$ \,and\, $T$ is a closed subset of \,$x$\, , then \,$A\,,\,B\,,\,S$ \,and\, \,$T$\, have a unique common fixed point in \,$x$.
\end{Theorem}

{\bf Proof.} Suppose that \,$(\,B\,,\,T\,)$\, satisfies property \,$(\,E.A.\,),$\, then there exists a
sequence \,$\{x_{\,n}\}$\, in \,$x$\, such that
\[\mathop {\lim }\limits_{n\,
\, \to \,\,\infty }Bx_{\,n}\,=\,\mathop {\lim }\limits_{n\,\, \to \,\,\infty }\,Tx_{\,n}\,=\,
z{\hspace{0.5cm}}\,for\;\;some \;\; z \,\in\, X.\]
Therefore , we have
\[\mathop {\lim }\limits_{n\,
\, \to \,\,\infty }\mu(\,Bx_{\,n} \,,\, Tx_{\,n} \,,\, t\,) \,=\, 1 \hspace{2.5cm}\]
Since \,$B(\,x\,)\,\subset\,S(\,x\,)$ , there exists in \,$X$\, a sequence \,$\{y_{\,n}\}$\, such that
\\$Bx_{\,n}\,=\,Sy_{\,n}$ . Putting \,$x\,=\,y_{\,n}$\, and \,$y\,=\,x_{\,n}$ in $(\,4\,)$
\\\\${\hspace{0.5cm}} \,F(\,\mu(\,Ay_{\,n}\,,\,Bx_{\,n}\,,\,t\,)\,,\,\mu(\,Sy_{\,n}\,,\,Tx_{\,n}\,,\,t\,)\,,\,
\mu(\,Ay_{\,n}\,,\,Sy_{\,n}\,,\,t\,) \,,\,$
\[{\hspace{3.5cm}}\mu(\,Bx_{\,n}\,,\,Tx_{\,n}\,,\,t\,)\,,\,1\,,\,
\mu(\,Ay_{\,n}\,,\,Tx_{\,n}\,,\,t\,)\,) \;\leq\; 1\,\,
{\hspace{0.5cm}}\forall\; t \,>\, 0\]
\[\Longrightarrow \;\; F(\,\mu(\,Ay_{\,n}\,,\,Bx_{\,n}\,,\,t\,)\,,\,
\mu(\,Bx_{\,n}\,,\,Tx_{\,n}\,,\,t\,)\,,\,
\mu(\,Ay_{\,n}\,,\,Bx_{\,n}\,,\,t\,)\,,
{\hspace{2.5cm}}\]
\[{\hspace{5.5cm}}\mu(\,Bx_{\,n}\,,\,Tx_{\,n}\,,\,t\,)\,,\,1\,,\,
\mu(\,Ay_{\,n}\,,\,Tx_{\,n}\,,\,t\,)\,) \;\leq\; 1\,\,\]
\[\Longrightarrow\;\; F(\,\mu(\,Ay_{\,n}\,,\,Bx_{\,n}\,,\,t\,)\,,\,\mu(\,Bx_{\,n}
\,,\,Bx_{\,n}\,,\,t\,)\,,\,
\mu(\,Ay_{\,n}\,,\,Bx_{\,n}\,,\,t\,)\,,
{\hspace{3.0cm}}\]
\[{\hspace{5.5cm}}\mu(\,Bx_{\,n}\,,\,Bx_{\,n}\,,\,t\,)\,,\,1\,,\,
\mu(\,Ay_{\,n}\,,\,Bx_{\,n}\,,\,t\,)\,) \;\leq\; 1\,\,\]
Taking the limit as \,$n \,\longrightarrow\, \infty, $\, we obtain :
\[{\hspace{0.6cm}}F(\,\mu(\,Ay_{\,n}\,,\,Bx_{\,n}\,,\,t\,)\,,\,1\,,\,
\mu(\,Ay_{\,n}\,,\,Bx_{\,n}\,,\,t\,)\,,\,
1\,,\,1\,,\,\mu(\,Ay_{\,n}\,,\,Bx_{\,n}\,,\,t\,)\,) \;\leq\; 1\,\,\]
which is a contradiction of \,$(\,F_{\,2}\,)$ , then we have
\[\mathop {\lim }\limits_{n\,
\, \to \,\,\infty }\mu(\,Ay_{\,n}\,,\,Bx_{\,n}\,,\,t\,) \,=\, 1\]
\[\Longrightarrow\;\; \mathop {\lim }\limits_{n\,
\, \to \,\,\infty }\,Ay_{\,n} \,=\, z \hspace{3.2cm}\]
Suppose that \,$S(\,x\,)$\, is a closed subspace of \,$X$ .
Then \,$z \,=\, Su$\, for some \,$u \,\in\, x.$\,
 Putting \,$x \,=\, u\,$\, and \,$y \,=\, x_{\,n}$\, in \,$(\,4\,)$\, we obtain :
\[ \,F(\mu(\,Au\,,\,Bx_{\,n}\,,\,t\,)\,,\,\mu(\,Su\,,\,Tx_{\,n}\,,\,t\,)\,,\,
\mu(\,Au\,,\,Su\,,\,t\,)\,,\,
\mu(\,Bx_{\,n}\,,\,Tx_{\,n}\,,\,t\,)\,,\]
\[{\hspace{6.5cm}}\mu(\,Su\,,\,Bx_{\,n}\,,\,t\,)\,,\,
\mu(\,Au\,,\,Tx_{\,n}\,,\,t\,)) \;\leq\; 1\,\]
Letting \,$n \,\longrightarrow\, \infty$\, we have
\[ F(\,\mu(\,Au\,,\,z\,,\,t\,)\,,\,1\,,\,\mu(\,Au\,,\,z\,,\,t\,)\,,\,1\,,\,1\,,\,
\mu(\,Au\,,\,z\,,\,t\,)\,) \;\leq\; 1\,\]
which is a contradiction of \,$(\,F_{\,2}\,).$\, Hence ,
\[ \mu(\,Au\,,\,z\,,\,t\,) \,=\, 1\]
\[\Longrightarrow\;\; Au\,=\,z \,\;\Longrightarrow\; \, z\,=\, Au \,=\, Su{\hspace{0.5cm}}\]
Since \,$A(\,X\,) \,\subset\, T(\,X\,)\,,$\, there exists \,$v \,\in \, x$\, such that
$z \,=\, Au \,=\, Tv\,.$
If \,$Az \,\neq\, z$\, and putting \,$x \,=\, u$\, and \,$y \,=\, v$\, in \,$(\,4\,)$\,,
 then we get
\[F\,(\,\mu(\,Au\,,\,Bv\,,\,t\,)\,,\,\mu(\,Su\,,\,Tv\,,\,t\,)\,,\,
\mu(\,Au\,,\,Su\,,\,t\,)\,,\,\mu(\,Bv\,,\,Tv\,,\,t\,),\]
\[{\hspace{6.5cm}}\mu(\,Su\,,\,Bv\,,\,t\,)\,,\,\mu(\,Au\,,\,Tv\,,\,t\,)\,) \;\leq\; 1\,\]
\[F(\,\mu(\,z\,,\,Bv\,,\,t\,)\,,\,1\,,\,1\,,\,\mu(\,z\,,\,Bv\,,\,t\,)\,,\,
\mu(\,z\,,\,Bv\,,\,t\,)\,,\,1\,) \;\leq\; 1\,{\hspace{1.3cm}}\]
a contradiction of \,$(\,F_{\,1}\,)$\, then
\[ Bv \,=\, z \, \;\Longrightarrow\;\, Tv \,=\, Bv \,=\, z\]
\[\Longrightarrow\;\; Au \,=\, Su \,=\, z \,=\, Tv \,=\, Bv \hspace{1.5cm}\]
Since the pair \,$(\,A\,,\,S\,)$\, is weakly compatible , we have
\[\, ASu \,=\, SAu \;\;\Longrightarrow\;\; Az \,=\, Sz\]
If \,$Az \,\neq\, z$\, and putting \,$x \,=\, z \,=\, y$\, in \,$(\,4\,)$\,
\[F(\,\mu(\,Az\,,\,Bv\,,\,t\,)\,,\,\mu(\,Sz\,,\,Tv\,,\,t\,)\,,\,
\mu(\,Az\,,\,Sz\,,\,t\,)\,,\,\mu(\,Bv\,,\,Tv\,,\,t\,),\]
\[{\hspace{6.5cm}}\mu(\,Sz\,,\,Bv\,,t\,)\,,\,\mu(\,Az\,,\,Tv\,,\,t\,)\,) \;\leq\; 1 \]
\[{\hspace{0.5cm}}F(\,\mu(\,Az\,,\,z\,,\,t\,)\,,\,\mu(\,Az\,,\,z\,,\,t\,)\,,\,1\,,\,1\,,\,
\mu(\,Az\,,\,z\,,\,t\,)\,,\,\mu(\,Az\,,\,z\,,\,t\,)\,) \;\leq\; 1 \]
which is a contradiction of \,$(\,F_{\,3\,}\,)$\,. Then \,$Az \,=\, Sz \,=\, z$\,
\\Since the pair \,$(\,B\,,\,T\,)$\, is weakly compatible , we have
\[BTv \,=\, TBv \hspace{0.5cm}i.e ,\hspace{0.5cm} Bz \,=\, Tz\]
If \,$Bz \,\neq\, z$\, and putting \,$x \,=\, z \,=\, y$\, in $(\,4\,)$
\[F(\,\mu(\,Az\,,\,Bz\,,\,t\,)\,,\,\mu(\,Sz\,,\,Tz\,,\,t\,)\,,\,
\mu(\,Az\,,\,Sz\,,\,t\,)\,,\,\mu(\,Bz\,,\,Tz\,,\,t\,),\]
\[{\hspace{6.5cm}}\mu(\,Sz\,,\,Bz\,,\,t\,)\,,\,\mu(\,Az\,,\,Tz\,,\,t\,)\,) \;\leq\; 1\]
\[{\hspace{0.5cm}} F(\,\mu(\,z\,,\,Bz\,,\,t\,)\,,\,\mu(\,z\,,\,Bz\,,\,t\,)\,,\,1\,,\,1\,,\,
\mu(\,z\,,\,Bz\,,\,t\,)\,,\,\mu(\,z\,,\,Bz\,,\,t\,)\,) \;\leq\; 1\]
which is a contradiction of \,$(\,F_{\,3}\,)$.
\\Hence $z \,= \,Bz\, = \,Tz \,=\, Az\, =\, Sz$ and \,$z$\, is common fixed point of
$A \,,\, B \,,\, S$ \,and\, $T.$
\\Suppose that $A\,,\,B\,,\,S$ \,and\, $T$\, have another fixed point $z_{\,1}.$
\\Then , we put $x=z$\, and \,$y\,=\,z_{\,1}$\, in \,$(\,4\,)$
\[F(\,\mu(\,Az\,,\,Bz_{\,1},t\,)\,,\,\mu(\,Sz\,,\,Tz_{\,1}\,,\,t\,)\,,\,
\mu(\,Az\,,\,Sz\,,\,t\,),
\mu(\,Bz_{\,1}\,,\,Tz_{\,1}\,,\,t\,)\,,\]
\[{\hspace{6.5cm}}\mu(\,Sz\,,\,Bz_{\,1}\,,\,t\,)\,,\,
\mu(\,Az\,,\,Tz_{\,1}\,,\,t\,)\,) \;\leq\; 1 \]
\[F(\mu(\,z\,,\,z_{\,1}\,,\,t\,)\,,\,\mu(\,z\,,\,z_{\,1}\,,\,t\,)\,,\,1\,,\,1\,,\,
\mu(\,z\,,\,z_{\,1}\,,\,t\,)\,,\,\mu(\,z\,,\,z_{\,1}\,,\,t\,)\,) \;\leq\; 1 \hspace{0.4cm}\]
which is a contradiction of \,$(\,F_{\,3}\,)$ , then we get \,$z \,=\, z_{\,1}$ .

\begin{Corollary}Let \,$A \,,\, B\,,\, S$\, and \,$T$\, be self-mappings of a fuzzy metric space  \,$(\,X \,,\, \mu \,,\, \ast\,)$\, satisfying the following conditions.
\[F(\,\mu(\,Ax\,,\,By\,,\,t\,)\,,\,\mu(\,Sx\,,\,Ty\,,\,t\,)\,,\,
\mu(\,Ax\,,\,Sx\,,\,t\,)\,,\,\mu(\,By\,,\,Ty\,,\,t\,),\]
\[{\hspace{6.5cm}}\mu(\,Sx\,,\,By\,,\,t\,)\,,\,\mu(\,Ax\,,\,Ty\,,\,t\,) \;\leq\; 1 \]
for all \,$x \,,\, y \;in\; X$ and where \,$F$\, satisfies property $(\,F_{\,1}\,) \,,\, (\,F_{\,2}\,)$ \,and\, $(\,F_{\,3}\,).$
\\Suppose that \,$(\,A \,,\, S\,)$\, or \,$(\,B \,,\, T\,)$\, satisfies property \,$(\,E.A\,)$\, and the pairs \,$(\,A \,,\, S\,)$\, and \,$(\,B \,,\, T\,)$\, are weakly compatible. If \,$S(\,X\,)$\, and \,$T(\,X\,)$\, are closed subset of \,$X$\,, then \,$A \,,\, B\,,\, S$\, and \,$T$\, have a unique common fixed point in \,$X$.
\end{Corollary}

\section{Well-Posedness Of Common Fixed Point \\Theorem}


\begin{Definition}
Let \,$(\,X \,,\, \mu \,,\, \ast\,)$\, be fuzzy metric space and \,$\mathcal{P}$\, a set of self-mappings of \,$X$. The common fixed point problem of the set \,$\mathcal{P}$\, is said to be well-posed if
\\$ {\hspace{1.0cm}}(\,1\,)$ \,$\mathcal{P}$\, has a unique common fixed point \,$x$\, in \,$X$. ( that is , \,$x$\, is the unique point in \,$X$\, such that $Ax \,=\, x $\, \,$\forall\;A \,\in\,\mathcal{P}$
\\$ {\hspace{1.0cm}}(\,2\,)$ For every sequence \,$\{x_{\,n}\}$\, in \,$X$\, such that
\[\mathop {\lim }\limits_{\,n\,\, \to \,\,\infty }\,\mu(\,x_{\,n} \,,\, Ax_{\,n} \,,\, t\,) \,=\, 1 \;\;\forall\,\; A \,\in\, \mathcal{P}\]
\[\Longrightarrow\;\; \mathop {\lim }\limits_{\,n\,\, \to \,\,\infty }\,\mu(\,x_{\,n} \,,\, x \,,\, t\,) \,=\, 1 \hspace{3.5cm}\]
\end{Definition}

\begin{Theorem}
Let \,$A\,,\,B$\, and \,$I$\, be three self mappings of a fuzzy metric space
\,$(\,X \,,\, \mu \,,\, \ast\,)$\, such that :
\\${\hspace{2.0cm}}(\,i\,)$ The pair \,$\{\,A\,,\,I\,\}$\, and \,$\{B\,,\,I\}$\, are weakly compatible,
\\${\hspace{2.0cm}}(\,ii\,)$ The mappings \,$A\,,\,B$\, and \,$I$\, satisfy the property $(\,E.A.\,)\,,$
\\${\hspace{2.0cm}}(\,iii\,) \,F\left(\,\mu\left(\,Ax\,,\,By\,,\,\frac{t}{2}\,\right)\,,\,
\mu\left(\,Ix\,,\,Iy\,,\,\frac{t}{2}\,\right)\,,\,
\mu\left(\,Ix\,,\,Ax\,,\,\frac{t}{2}\,\right)\right.\,,\,$
\\${\hspace{2.0cm}}\left.\mu\left(\,Iy\,,\,By\,,\,\frac{t}{2}\,\right)\,,\,
\mu\left(\,Ix\,,\,By\,,\,\frac{t}{2}\,\right)\,,\,
\mu\left(\,Iy\,,\,Ax\,,\,\frac{t}{2}\,\right)\right)\;
\leq\; 1\,$
\\$\forall\; x\;\neq\;y\; \; in \;\;X\,.$
Where \,$F$\, satisfies property $(\,F_{1}\,)\,,\,(\,F_{\,2}\,)\,,\,(\,F_{\,3}\,)\,and\,(\,F_{\,4}\,)$
\\${\hspace{2.0cm}}(\,iv\,)$ \,$I(\,x\,)$\,is closed ,
\\Then the common fixed point problem of \,$A\,,\,B$\, and \,$I$\, is well posed.
\end{Theorem}

{\bf Proof.} Let \,$\{x_{\,n}\}$\, be a sequence in \,$X$\, such that
\[\mathop {\lim }\limits_{n\,\, \to \,\,\infty }\,\mu\left(\,x_{\,n} \,,\, Ax_{\,n} \,,\, \frac{t}{2}\,\right) \,=\,
\mathop {\lim }\limits_{n\,\, \to \,\,\infty }\mu\left(\,x_{\,n} \,,\, Bx_{\,n} \,,\, \frac{t}{2}\,\right) \,=\,
\mathop {\lim }\limits_{n\,\, \to \,\,\infty }\,\mu\left(\,x_{\,n} \,,\, Ix_{\,n} \,,\, \frac{t}{2}\,\right) \,=\,1 \]
Putting \,$y\,=\,x_{\,n}$\, in \,$(\,iii\,)$\, then
\[F\left(\mu\left(\,Ax \,,\, Bx_{\,n} \,,\, \frac{t}{2}\,\right) \,,\,
\mu\left(\,Ix \,,\, Ix_{\,n} \,,\, \frac{t}{2}\,\right) \,,\,
\mu\left(\,Ix \,,\, Ax \,,\, \frac{t}{2}\,\right)\,
\right.,\]
\[{\hspace{2.0cm}}\left.\mu\left(\,Ix_{\,n} \,,\, Bx_{\,n} \,,\, \frac{t}{2}\,\right)\,,\,
\mu\left(\,Ix \,,\, Bx_{\,n} \,,\, \frac{t}{2}\,\right) \,,\,
\mu\left(\,Ix_{\,n} \,,\, Ax \,,\, \frac{t}{2}\,\right)\right) \;\leq\; 1\]
\[\Longrightarrow\;\; F\left(\,\mu\left(\,x \,,\, Bx_{\,n} \,,\, \frac{t}{2}\,\right) \,,\,
\mu\left(\,x \,,\, Ix_{\,n} \,,\, \frac{t}{2}\,\right) \,,\,
\mu\left(\,x \,,\, x\,,\, \frac{t}{2}\,\right)\,
\right.,{\hspace{3.5cm}}\]
\[{\hspace{2.4cm}}\left.\mu\left(\,Ix_{\,n} \,,\, Bx_{\,n} \,,\, \frac{t}{2}\,\right) \,,\ \mu\left(\,x \,,\, Bx_{\,n} \,,\, \frac{t}{2}\,\right) \,,\,
\mu\left(\,Ix_{\,n} \,,\, x \,,\, \frac{t}{2}\,\right)\,\right) \;\leq\; 1\]
\[\Longrightarrow\;\; F\left(\,\mu\left(\,x \,,\, Bx_{\,n} \,,\, \frac{t}{2}\,\right) \,,\, \mu\left(\,x \,,\, Ix_{\,n} \,,\, \frac{t}{2}\,\right)\right. \,,\, 1 \,,\,
\mu\left(\,Ix_{\,n} \,,\, Bx_{\,n} \,,\, \frac{t}{2}\,\right) \,,\,{\hspace{2.0cm}}\]
\[{\hspace{6.4cm}}\left.\mu\left(\,x \,,\, Bx_{\,n} \,,\, \frac{t}{2}\,\right) \,,\,
\mu\left(\,Ix_{\,n} \,,\, x \,,\, \frac{t}{2}\,\right)\right) \;\leq\; 1\]
By \,$(\,F_{\,4\,}\,)$\, we get
\[\mu\left(\,x \,,\, Bx_{\,n} \,,\, \frac{t}{2}\,\right) \;\;\geq\;\; \frac{1}{k}\;
\mu\left(\,Ix_{\,n} \,,\, Bx_{\,n} \,,\, \frac{t}{2}\,\right)\]
Therefore ,
\[\mu\left(\,x \,,\, x_{\,n} \,,\, t\,\right) \;\;\geq\;\;
\mu\left(\,x \,,\, Bx_{\,n} \,,\, \frac{t}{2}\,\right) \,\ast\,
\mu\left(\,Bx_{\,n} \,,\, x_{\,n} \,,\, \frac{t}{2}\,\right){\hspace{3.2cm}}\]
\[\hspace{0.8cm}\geq \;\; \frac{1}{k}\;\left\{\mu\left(\,Ix_{\,n} \,,\, Bx_{\,n} \,,\,
\frac{t}{2}\,\right)\right\} \,\ast\,
\mu\left(\,Bx_{\,n} \,,\, x_{\,n} \,,\, \frac{t}{2}\,\right)\]
\[{\hspace{2.8cm}}\geq\;\; \frac{1}{k}\;\left\{\,\mu\left(\,Ix_{\,n} \,,\, x_{\,n} \,,\,
\frac{t}{4}\,\right) \,\ast\, \mu\left(\,Bx_{\,n} \,,\, x_{\,n} \,,\,
\frac{t}{4}\,\right)\right\}\,
\ast\, \mu\left(\,Bx_{\,n} \,,\,x_{\,n} \,,\,\frac{t}{2}\,\right)\]
taking limit as \,$n\,\rightarrow\,\infty$ , we have
\[\mathop {\lim }\limits_{n\,\, \to \,\,\infty }\mu(\,x \,,\, x_{\,n} \,,\, t\,) \;\geq\; \frac{1}{k} \;>\; 1{\hspace{3.2cm}}\]
\[\mathop {\lim }\limits_{n\,\, \to \,\,\infty }\mu(\,x \,,\, x_{\,n} \,,\, t\,) \;=\; 1{\hspace{4.2cm}}\]
This completes the proof.

\begin{Theorem}Let \,$A\,,\,B\,,\,S $\,and \,$T$\, be self - mappings of a fuzzy metric space \,$(\,X \,,\, \mu \,,\, \ast\,)$\, satisfying the following conditions.
\[A(\,X\,) \,\subset\, T(\,X\,) \hspace{0.5cm}and\hspace{0.5cm} B(\,X\,) \,\subset\, S(\,X\,)\]
\[{\hspace{0.5cm}} \,F\left(\,\mu\left(\,Ax \,,\, By \,,\, \frac{t}{2}\,\right) \,,\,
\mu\left(\,Sx \,,\, Ty \,,\, \frac{t}{2}\,\right) \,,\,
\mu\left(\,Ax \,,\, Sx \,,\, \frac{t}{2}\,\right) \,,\,
\mu\left(\,By \,,\, Ty \,,\, \frac{t}{2}\,\right)\right.,\]
\[{\hspace{4.9cm}}\left.\mu\left(\,Sx \,,\, By \,,\, \frac{t}{2}\,\right) \,,\,
\mu\left(\,Ax \,,\, Ty \,,\, \frac{t}{2}\,\right)\right) \;\leq\; 1 {\hspace{1.0cm}}\cdots\hspace{0.5cm}(\,5\,)\]
\,for all \,$x \,,\, y \,\in\,X$ and where \,$F$\, satisfies property $(\,F_{\,1}\,) \,,\, (\,F_{\,2}\,) \,,\, (\,F_{3}\,)$\,and\,$(\,F_{\,4}\,).$
\\Suppose that $(\,A \,,\, S\,)$ \,or\, $(\,B \,,\, T\,)$\, satisfies property $(\,E.A.\,)$ and the pairs $(\,A \,,\, S\,)$ \,and\, $(\,B \,,\, T\,)$ are weakly compatible. If the range of one $A \,,\, B\,,\, S$\,and \,$T$\, is a closed subset of \,$X$\, , then the common fixed point problem of $A \,,\,B \,,\, S$ \,and \,$T$\, are well posed .
\end{Theorem}

{\bf Proof.} Suppose that $(\,B \,,\, T\,)$ satisfies property \,$(\,E.A\,)$\,, then there exists a sequence \,$\{\,x_{\,n}\,\}$\, in \,$X$\, such that
\[\mathop {\lim }\limits_{n\,\, \to \,\,\infty }\mu\left(\,x_{\,n}\,,\,Bx_{\,n}\,,\,\frac{t}{2}\,\right) \,=\, \mathop {\lim }\limits_{n\,\, \to \,\,\infty }\mu\left(\,x_{\,n}\,,\,Tx_{\,n}\,,\,\frac{t}{2}\,\right)\,=\,1\]
putting \,$y\,=\,x_{\,n}$\, in \,$(\,5\,)$\,, we have
\[F\left(\,\mu\left(\,Ax\,,\,Bx_{\,n}\,,\,t\,\right)\,,\,
\mu\left(\,Sx\,,\,Tx_{\,n}\,,\,\frac{t}{2}\,\right),
\mu\left(\,Ax\,,\,Sx\,,\,\frac{t}{2}\,\right)\,,\,\right. \hspace{2.5cm}\]
\[{\hspace{1.5cm}}\left.\mu\left(\,Bx_{\,n}\,,\,Tx_{\,n}\,,\,\frac{t}{2}\,\right)\,,\,
\mu\left(\,Sx\,,\,Bx_{\,n}\,,\,\frac{t}{2}\,\right)\,,\,
\mu\left(\,Ax\,,\,Tx_{\,n}\,,\,\frac{t}{2}\,\right)\right) \;\leq\; 1 \]
\[F\left(\,\mu\left(\,x\,,\,Bx_{\,n}\,,\,\frac{t}{2}\,\right)\,,\,
\mu\left(\,x\,,\,Tx_{n}\,,\,\frac{t}{2}\,\right)\,,\,1\,,\,
\mu\left(\,Bx_{\,n}\,,\,Tx_{\,n}\,,\,\frac{t}{2}\,\right)\right. \hspace{2.5cm}\]
\[{\hspace{6.1cm}}\left.\mu\left(\,x\,,\,Bx_{\,n}\,,\,\frac{t}{2}\,\right)\,,\,
\mu\left(\,x\,,\,Tx_{\,n}\,,\,\frac{t}{2}\,\right)\right) \;
\leq\; 1 \]
By $(\,F_{\,4}\,)$ we get
\[\mu\left(\,x\,,\,Bx_{\,n}\,,\,\frac{t}{2}\,\right) \;\geq\; \frac{1}{k}\;
\mu\left(\,Bx_{\,n}\,,\,Tx_{\,n}\,,\,\frac{t}{2}\,\right) .\]
Therefore ,
\[\mu\left(\,x\,,\,x_{\,n}\,,\,t\,\right) \;\geq\;
\mu\left(\,x\,,\,Bx_{\,n}\,,\,\frac{t}{2}\,\right) \,\ast\,
\mu\left(\,Bx_{\,n}\,,\,x_{\,n}\,,\,\frac{t}{2}\,\right){\hspace{3.5cm}}\]
\[\hspace{0.5cm}\geq\; \frac{1}{k}\;\left\{\,\mu\left(\,Bx_{\,n}\,,\,Tx_{\,n}\,,\,
\frac{t}{2}\,\right)\right\} \,\ast\,
\mu\left(\,Bx_{\,n}\,,\,x_{\,n}\,,\,\frac{t}{2}\,\right)\]
\[{\hspace{2.7cm}}\geq\; \frac{1}{k}\;\left\{\,\mu\left(\,Bx_{\,n}\,,\,x_{\,n}\,,\,
\frac{t}{4}\,\right)\,
\ast\, \mu\left(\,Tx_{\,n}\,,\,x_{\,n}\,,\,\frac{t}{4}\,\right)\,\right\} \,\ast\,
\mu\left(\,Bx_{\,n}\,,\,x_{\,n}\,,\frac{t}{2}\,\right)\]
taking limit as \,$n\,\longrightarrow\,\infty$\,
\[\mathop {\lim }\limits_{n\,\, \to \,\,\infty }\mu(\,x\,,\,x_{\,n}\,,\,t\,) \;\geq\; \frac{1}{k} \,>\, 1{\hspace{3.4cm}}\]
\[\mathop {\lim }\limits_{n\,\, \to \,\,\infty }\mu(\,x\,,\,x_{\,n}\,,\,t\,) \,=\, 1{\hspace{4.2cm}}\]
\\This completes the proof.

\begin{Corollary}Let \,$A \,,\, B \,,\, S$\, and \,$T$\, be self-mappings of a fuzzy metric space  \,$(\,X \,,\, \mu \,,\, \ast\,)$\, satisfying the following conditions.
\[F\left(\,\mu\left(\,Ax\,,\,By\,,\,\frac{t}{2}\,\right)\,,\,
\mu\left(\,Sx\,,\,Ty\,,\,\frac{t}{2}\,\right)\,,\,
\mu\left(\,Ax\,,\,Sx\,,\,\frac{t}{2}\,\right)\,,\,
\mu\left(\,By\,,\,Ty\,,\,\frac{t}{2}\,\right)\right.,\]
\[{\hspace{5.5cm}}\left.\mu\left(\,Sx\,,\,By\,,\,\frac{t}{2}\,\right)\,,\,
\mu\left(\,Ax\,,\,Ty\,,\,\frac{t}{2}\,\right)\,\right) \;\leq\; 1 \]
for all \,$x\,,\,y\,\in\,X$ and where \,$F$\, satisfies property $(\,F_{\,1}\,) \,,\, (\,F_{\,2}\,) \,,\, (\,F_{\,3}\,)$ \,and\, $(\,F_{\,4}\,).$\,
Suppose that \,$(\,A\,,\,S\,)$\, or \,$(\,B\,,\,T\,)$\, satisfies property \,$(\,E.A\,)$\, and the pairs \,$(\,A \,,\, S\,)$\, and \,$(\,B \,,\, T\,)$\, are weakly compatible. If \,$S(\,x\,)$\, and \,$T(\,x\,)$\, are closed subset of \,$X$\,, then the common fixed point problem of \,$A \,,\, B\,,\, S$\, and \,$T$\, are well posed .
\end{Corollary}

\end{document}